\documentclass[number,citesort,MSNbibl,seceqn,dvips]{arxbj}

\aid{0}
\volume{18}
\issue{2}
\pubyear{2012}
\firstpage{434}
\lastpage{454}
\doi{10.3150/10-BEJ341}

\makeatletter
\renewcommand{\epsilon}{\varepsilon}
\newtheorem{theorem}{Theorem}
\newtheorem{lemma}[theorem]{Lemma}
\newtheorem{proposition}[theorem]{Proposition}
\newproclaim{definition}{Definition}
\makeatother

\begin{document}
\begin{frontmatter}

\title{Thermodynamics and concentration}
\runtitle{Thermodynamics and concentration}

\begin{aug}
\author{\fnms{Andreas} \snm{Maurer}\corref{}\ead[label=e1]{am@andreas-maurer.eu}}
\runauthor{A. Maurer}
\address{Adalbertstr. 55, D-80799 M\"{u}nchen, Germany. \printead{e1}}
\end{aug}

\received{\smonth{8} \syear{2009}}
\revised{\smonth{5} \syear{2010}}

\begin{abstract}
We show that the thermal subadditivity of entropy provides a common basis
to derive a strong form of the bounded difference inequality and related
results as well as more recent inequalities applicable to convex
Lipschitz functions, random symmetric matrices, shortest travelling
salesmen paths and weakly self-bounding functions. We also give two new
concentration inequalities.
\end{abstract}


\vspace*{-2pt}
\end{frontmatter}

\section{Introduction}\vspace*{-2pt}

Concentration inequalities bound the probabilities that random
quantities deviate from their average, median or otherwise typical
values. They are at the heart of empirical science and play an
important role in the study of natural and artificial learning systems.

An early concentration inequality for sums was given by Chebychev and
Bienaym\'{e} in the 19th century \cite{Chebychev1874} and allowed a
rigorous proof of the weak law of large numbers. The subject has since
been developed by Bernstein, Chernoff, Bennett, Hoeffding and many
others \cite{Bernstein1927,Hoeffding1963}, and results were extended
from sums to more general and complicated nonlinear functions. During
the past few decades, research activity has been stimulated by the
contributions of Michel Talagrand \cite{Talagrand1995,Talagrand1996}
and by the relevance of concentration phenomena to the rapidly growing
field of computer science. Some concentration inequalities, like the
well-known bounded difference inequality, have become standard tools in
the analysis of algorithms \cite{McDiarmid1998}. Nevertheless, a
unified and elementary basis for the derivation of the many available
results is still missing.

One of the more recent methods used to derive concentration
inequalities, the so-called \textit{entropy method}, is rooted in the
early investigations of Boltzmann \cite{Boltzmann1877} and Gibbs
\cite{Gibbs1902} into the foundations of statistical mechanics. A
general problem of statistical mechanics is to demonstrate the
``equivalence of ensembles'', which can be interpreted as an exponential
concentration property of the Hamiltonian, or energy function. While
the modern entropy method evolved along a complicated historical path
via quantum field theory and the logarithmic Sobolev inequality of
Leonard Gross \cite{Gross1975}, its hidden simplicity was understood
and emphasized by Michel Ledoux, who also recognized the key role that
the subadditivity of entropy can play in the derivation of
concentration inequalities \cite{Ledoux1996,Ledoux2001}. Recently,
Boucheron \textit{et al.} \cite{Boucheron2009} showed that the entropy method
is sufficiently strong to derive a~form of Talagrand's convex distance
inequality.\vadjust{\goodbreak}

The purpose of this paper is to advertise the subadditivity of entropy
as a unified basis for the derivation of concentration inequalities for
functions on product spaces and to demonstrate the benefits of
formulating the concentration problem in the language of statistical
thermodynamics, an approach proposed by David McAllester
\cite{McAllester2002}.

Our method consists of three steps. The first step (Theorem
\ref{Theorem basic Herbst argument}) expresses the log-Laplace
transform (or, more directly, the deviation probability) in terms of an
integral of the thermal entropy over a range of inverse temperatures.
This step encapsulates the so-called Herbst argument.

The second step (Theorem \ref{Theorem Entropy subadditivity}) is the
tensorization inequality, or, more properly, a~thermal subadditivity
property of entropy. It asserts that the entropy of a~system is no
greater than the thermal average of the sum of entropies of the
constituent subsystems.

The third step (Theorem \ref{Theorem Entropy fluctuation integral})
expresses the entropy of the subsystem in terms of thermal energy
fluctuations.

All three steps are elementary and their combination leads to a general
concentration result (Theorem \ref{Theorem general concentration}) that
can be used whenever we succeed in controlling the latter
fluctuations.

We then use the method to first derive a strong form of the bounded
difference inequality and an inequality given by McDiarmid and related
to Bennett's inequality \cite{McDiarmid1998}. These results are
normally not associated with the entropy method. Then monotonicity
properties of thermal energy fluctuations, or bounds thereof, are
exploited to derive two apparantly novel sub-Gaussian tail-bounds and
to give a new proof of an upper tail-bound in \cite{Maurer2006} that
improves on some results obtained from Talagrand's convex distance
inequality. Finally, we show how our method can be extended in a
generic way using self-boundedness and/or decoupling, and illustrate
this extension by deriving a concentration inequality that underlies
the recent new proof of the convex distance inequality~\cite{Boucheron2009}.

Clearly statement and proof of all the results presented in this paper
would be possible on a~purely formal basis without any reference to
physics. The author believes, however, that positioning the subject in
a broader scientific context highlights its historical connections and
gives access to a valuable source of intuition.

In the next section, we describe the connection between entropy and
concentration and introduce several thermodynamic functions. We then
transfer these concepts to product spaces and present the tensorization
inequality. The remaining sections are dedicated to applications, and
we conclude with a tabular summary of the notation used in the paper.

\section{Entropy and concentration}

Let $( \Omega ,\Sigma ,\mu ) $ be a probability space and $%
f\in L_{\infty }[ \mu ] $ be a fixed function whose
concentration properties are to be studied.

We interpret the points $x\in \Omega $ as possible states of a physical
system and $f$ as the negative energy (or Hamiltonian) function, so that $%
-f( \mathbf{x}) $ is the system's energy in the state $x$.
The measure $\mu $ models an a priori probability distribution of
states in the absence of any constraining information.

We will ignore questions of measurability. If it seems necessary to the
reader, $\Omega $ may be taken as a potentially very large finite set,
the cardinality of which will play no role in our results. The
boundedness assumption is a simplification that is justified by the
fact that most of our results are vacuous for $f\notin L_{\infty
}[ \mu ] $. In the remaining cases, we will mention optimal
conditions on $f$.

For any $g\in L_{\infty }[ \mu ] $ we write $E[
g]
=\int_{\Omega }g\,\mathrm{d}\mu $ and $\sigma ^{2}[ g] =E[ ( g-E%
[ g] ) ^{2}] $.

\subsection{Thermal equilibrium and thermodynamic functions}

Our function $f$ defines a one-parameter family $\{ E_{\beta
f}\dvt\beta \in
\mathbb{R}
\} $ of expectation functionals by
\[
E_{\beta f}[ g] =\frac{E[ g\mathrm{e}^{\beta f}] }{E[
\mathrm{e}^{\beta f}] },\qquad g\in L_{\infty }[ \mu ].
\]
In statistical thermodynamics, $E_{\beta f}[ g] $ is the \textit{%
thermal expectation} of the observable $g$ at temperature $T=1/\beta $.
The normalizing expectation is called the \textit{partition function},
\[
Z_{\beta f}=E[ \mathrm{e}^{\beta f}] .
\]
The corresponding probability measure on $\Omega $,
\[
\mathrm{d}\mu _{\beta f}=Z_{\beta f}^{-1}\mathrm{e}^{\beta f}\,\mathrm{d}\mu ,
\]%
is called the canonical ensemble. It describes a system in thermal
equilibrium with a heat reservoir at temperature $T=1/\beta $. The
canonical ensemble has the density $\rho =Z_{\beta f}^{-1}\mathrm{e}^{\beta f}$,
which maximizes the Kullback--Leibler divergence or relative entropy
$KL( \rho\,\mathrm{d}\mu,\mathrm{d}\mu ) :=E[ \rho \ln \rho ] $,
given the expected internal energy $-E[ \rho f] $. The
parameter $\beta $ is the Lagrange multiplier corresponding to this
constraint. For a constant $c$ we have the obvious and important
identity $E_{\beta ( f+c) }[ g] =E_{\beta
f}[ g] $.

The corresponding maximal value of the Kullback--Leibler divergence is
the
\textit{canonical entropy}
%
\begin{equation}\label{Entropy_definition}
S_{f}( \beta ) =KL( Z_{\beta f}^{-1}\mathrm{e}^{\beta f}\,\mathrm{d}\mu,\mathrm{d}\mu )
=\beta E_{\beta f}[ f] -\ln Z_{\beta f}.
\end{equation}
Note that $S_{-f}( \beta ) =S_{f}( -\beta ) $, a
simple but very useful fact to pass from upper to lower tails.

For $\beta \neq 0$ the Helmholtz free energy is defined by
\[
A_{f}( \beta ) =\frac{1}{\beta }\ln Z_{\beta f}.
\]
Dividing (\ref{Entropy_definition}) by $\beta $ and writing $U=E_{\beta f}%
[ f] $, we obtain the classical thermodynamic relation
\[
A=U-TS,
\]
which describes the macroscopically available energy $A$ as the
difference between the total expected energy $U$ and an energy portion
$TS$, which is inaccessible due to ignorance of the microscopic state.

By L'H\^{o}pital's rule, we have $\lim_{\beta \rightarrow 0}A_{f}(
\beta ) =E[ f] $, so the free energy $A_{f}$ extends
continuously
to $%
\mathbb{R}
$ by setting $A_{f}( 0) =E[ f] $. We find
\[
A_{f}^{\prime }( \beta ) =\frac{1}{\beta }E_{\beta f}[ f%
] -\frac{1}{\beta ^{2}}\ln Z_{\beta f}=\beta ^{-2}S_{f}(
\beta ) .
\]%
Integrating this identity from zero to $\beta $ and multiplying with
$\beta, $ we obtain:

\begin{theorem}\label{Theorem basic Herbst argument}
For any $\beta >0$ we have
\[
\ln E\bigl[ \mathrm{e}^{\beta ( f-Ef) }\bigr] =\beta \int_{0}^{\beta }%
\frac{S_{f}( \gamma ) }{\gamma ^{2}}\,\mathrm{d}\gamma
\]
and, for $t\geq 0$,
\[
\Pr \{ f-Ef>t\} \leq \exp \biggl( \beta \int_{0}^{\beta }\frac{%
S_{f}( \gamma ) }{\gamma ^{2}}\,\mathrm{d}\gamma -\beta t\biggr) .
\]
\end{theorem}

\begin{pf}
\begin{eqnarray*}
\ln E\bigl[ \mathrm{e}^{\beta ( f-Ef) }\bigr] &=&\ln Z_{\beta f}-\beta E%
[ f] =\beta \bigl( A_{f}( \beta ) -A_{f}(
0) \bigr) \\
&=&\beta \int_{0}^{\beta }A_{f}^{\prime }( \gamma )\,\mathrm{d}\gamma
=\beta \int_{0}^{\beta }\frac{S_{f}( \gamma ) }{\gamma ^{2}}\,\mathrm{d}\gamma .
\end{eqnarray*}
Combining this with Markov's inequality gives the second assertion.
\end{pf}

The theorem shows how bounds on the canonical entropy can lead to
concentration results. In the following we present ways to arrive at
such bounds.

\subsection{Entropy and energy fluctuations}

The \textit{thermal variance} of a function $g\in L_{\infty }[ \mu %
] $ is denoted $\sigma _{\beta f}^{2}( g) $ and defined by
\[
\sigma _{\beta f}^{2}( g) =E_{\beta f}\bigl[ ( g-E_{\beta f}%
[ g] ) ^{2}\bigr] =E_{\beta f}[ g^{2}]
-( E_{\beta f}[ g] ) ^{2}.
\]
For constant $c$ we have $\sigma _{\beta ( f+c) }^{2}[ g%
] =\sigma _{\beta f}^{2}[ g] $.

We first give some simple results pertaining to the derivatives of the
partition function and the thermal expectations.

\begin{lemma}\label{General Simple Lemma}
The following formulas hold:
\begin{longlist}
\item[1.] $\frac{\mathrm{d}}{\mathrm{d}\beta }( \ln Z_{\beta f}) =E_{\beta f}[ f%
] $.\vspace*{1pt}

\item[2.] If $h\dvtx\beta \mapsto h( \beta ) \in L_{\infty }[ \mu %
] $ is differentiable and $( \mathrm{d}/\mathrm{d}\beta ) h( \beta
) \in L_{\infty }[ \mu ] $, then
\[
\frac{\mathrm{d}}{\mathrm{d}\beta }E_{\beta f}[ h( \beta ) ] =E_{\beta f}%
[ h( \beta ) f] -E_{\beta f}[ h( \beta
) ] E_{\beta f}[ f] +E_{\beta f}\biggl[ \frac{\mathrm{d}}{%
\mathrm{d}\beta }h( \beta ) \biggr] .
\]

\item[3.] $\frac{\mathrm{d}}{\mathrm{d}\beta }E_{\beta f}[ f^{k}] =E_{\beta f}[
f^{k+1}] -E_{\beta f}[ f^{k}] E_{\beta f}[
f] .$

\item[4.] $\frac{\mathrm{d}^{2}}{\mathrm{d}\beta ^{2}}( \ln Z_{\beta f}) =\frac{\mathrm{d}}{\mathrm{d}\beta }%
E_{\beta f}[ f] =\sigma _{\beta f}^{2}[ f] $.
\end{longlist}
\end{lemma}

\begin{pf}
1 is immediate and 2 is a straightforward computation. 3 and 4 are
immediate consequences of 1 and 2.
\end{pf}

The thermal variance of $f$ itself corresponds to energy fluctuations.
The next theorem represents entropy as a double integral of such
fluctuations. The utility of this representation to derive
concentration results has been noted by David McAllester
\cite{McAllester2002}.

\begin{theorem}\label{Theorem Entropy fluctuation integral}
We have for $\beta >0$
\[
S_{f}( \beta ) =\int_{0}^{\beta }\int_{t}^{\beta }\sigma
_{sf}^{2} [ f]\,\mathrm{d}s\,\mathrm{d}t.
\]
\end{theorem}

\begin{pf}
Using the previous lemma and the fundamental theorem of calculus, we
obtain
the formulas
\[
\beta E_{\beta f}[ f] =\int_{0}^{\beta }E_{\beta f}[
f] \,\mathrm{d}t=\int_{0}^{\beta }\biggl( \int_{0}^{\beta }\sigma
_{sf}^{2}[ f]\,\mathrm{d}s+E[ f] \biggr)\,\mathrm{d}t
\]
and
\[
\ln Z_{\beta f}=\int_{0}^{\beta }E_{tf}[ f]\,\mathrm{d}t
=\int_{0}^{\beta }\biggl( \int_{0}^{t}\sigma _{sf}^{2}[ f]\,\mathrm{d}s+E[ f] \biggr)\,\mathrm{d}t,
\]
which we subtract to obtain
\begin{eqnarray*}
S_{f}( \beta ) &=&\beta E_{\beta f}[ f] -\ln
Z_{\beta f}=\int_{0}^{\beta }\biggl( \int_{0}^{\beta }\sigma
_{sf}^{2}[ f]\,\mathrm{d}s-\int_{0}^{t}\sigma _{sf}^{2}[ f]\,\mathrm{d}s\biggr)\,\mathrm{d}t \\
&=&\int_{0}^{\beta }\biggl( \int_{t}^{\beta }\sigma _{sf}^{2}[
f]\,\mathrm{d}s\biggr)\,\mathrm{d}t.
\end{eqnarray*}
\upqed\end{pf}

Since bounding\vspace*{2pt} $\sigma _{\beta f}^{2}[ f] $ is central to
our method, it is worth mentioning an interpretation in terms of heat
capacity, or specific heat. Recall that $-E_{\beta f}[
f] $ is the expected internal energy. The rate of change of this
quantity with temperature $T$ is the heat capacity. By conclusion 4 of
Lemma \ref{General
Simple Lemma} we have
\[
\frac{\mathrm{d}}{\mathrm{d}T}( -E_{\beta f}[ f] ) =\frac{1}{T^{2}}%
\sigma _{\beta f}^{2}[ f] ,
\]
which exhibits the proportionality of heat capacity and energy
fluctuations.

\subsection{A variational entropy bound}

While Theorem \ref{Theorem Entropy fluctuation integral} is just an
elementary way of rewriting the canonical entropy, the following lemma
is typically a strict inequality that leads to the modified
logarithmic\vadjust{\goodbreak}
Sobolev inequality proposed by Massart in \cite{Massart2000a}. To state
it, we define the real function
%
\begin{equation}\label{Definition of psi}
\psi ( t) =\mathrm{e}^{t}-t-1,
\end{equation}
which arises from deleting the first two terms in the power series
expansion of the exponential function.

\begin{lemma}
\label{Lemma variational entropy bound}If $c\in
\mathbb{R},
$ then
\[
S_{f}( \beta ) \leq E_{\beta f}\bigl[ \psi \bigl( -\beta
( f-c) \bigr) \bigr] .
\]
\end{lemma}

\begin{pf}
Using $\ln t\leq t-1,$ we get
\[
\beta f-\ln Z_{\beta f}=\beta ( f-c) +\ln \frac{\mathrm{e}^{\beta c}}{%
Z_{\beta f}}\leq \beta ( f-c) +\biggl( \frac{\mathrm{e}^{\beta
c}}{Z_{\beta f}}-1\biggr).
\]
Taking the thermal expectation then gives
\begin{eqnarray*}
S_{f}( \beta ) &\leq &E_{\beta f}[ \beta ( f-c) %
] +\frac{\mathrm{e}^{\beta c}}{Z_{\beta f}}-1 \\
&=&E_{\beta f}[ \beta ( f-c) ] +E\biggl[ \frac{%
\mathrm{e}^{-\beta ( f-c) }\mathrm{e}^{\beta f}}{Z_{\beta f}}\biggr] -1 \\
&=&E_{\beta f}\bigl[ \mathrm{e}^{-\beta ( f-c) }+\beta ( f-c) -1%
\bigr].
\end{eqnarray*}%
\upqed\end{pf}

\section{Product spaces}

We now assume that $\Omega =\prod_{k=1}^{n}\Omega _{k}$ and $\mu
=\bigotimes _{k=1}^{n}\mu _{k}$, where each $\mu _{k}$ is the probability
measure representing the distribution of some variable $X_{k}$ in the
space $\Omega
_{k}$, where all the~$X_{k}$ are assumed to be mutually independent. The $%
X_{k}$ are irrelevant for the derivation of our inequalities, but they
are convenient in the discussion of applications.

If $\mathbf{x}=( x_{1},\ldots,x_{n}) \in \Omega $ describes a
state of a physical system, we can think of $x_{k}\in \Omega _{k}$ as
the state of the $k$th subsystem, which may be a particle or a more
abstract object, such as a spin assigned to the vertex of a graph. The
a priori measure $\mu $ assigns independent probabilities $\mu _{k}$ to
the states of the
subsystems. If the total energy is a sum of energies of the subsystems, $%
f=\sum f_{k}$, with $f_{k}\in L_{\infty }[ \mu _{k}] $, then
this is also true for the canonical ensemble $Z_{\beta f}\mathrm{e}^{\beta
f}\,\mathrm{d}\mu $ corresponding to non-interaction of the subsystems.

\subsection{Conditional expectations}

For $\mathbf{x}\in \Omega $,  $1\leq k\leq n$ and $y\in \Omega _{k}$ we
use $\mathbf{x}_{y,k}$ to denote the vector in $\Omega $, which is
obtained by replacing $x_{k}$ with $y$. We also write, for $g\in
L_{\infty }[ \mu ] $,
\[
E_{k}[ g] ( \mathbf{x}) =\int_{\Omega
_{k}}g( \mathbf{x}_{y,k})\,\mathrm{d}\mu _{k}( y)
=\int_{\Omega _{k}}g( x_{1},\ldots,x_{k-1},y,x_{k+1},\ldots,x_{n})\,\mathrm{d}\mu _{k}( y) .
\]

The operator $E_{k}$ corresponds to an expectation conditional to all
variables with indices different to $k$. We denote with
$\mathcal{A}_{k}$ the sub-algebra of $L_{\infty }[ \mu ] $
consisting of those functions that are independent of the $k$th
variable. $E_{k}$ is evidently a linear projection onto
$\mathcal{A}_{k}$. Also, the $E_{k}$ commute amongst
each other and, for $h\in L_{\infty }[ \mu ] $ and $g\in \mathcal{%
A}_{k}$, we have
%
\begin{equation} \label{conditional_identity}
E[ [ E_{k}h] g] =E[ E_{k}[ hg]] =E[ hg] .
\end{equation}

Replacing the operator $E$ by $E_{k}$ leads to the definition of
conditional
thermodynamic quantities, all of which are now members of the algebra $%
\mathcal{A}_{k}$:

\begin{itemize}
\item the conditional partition function $Z_{k,\beta f}=E_{k}[
\mathrm{e}^{\beta f}] $,

\item the conditional thermal expectation\vspace*{2pt} $E_{k,\beta f}[ g]
=Z_{k,\beta f}^{-1}E_{k}[ g\mathrm{e}^{\beta f}] $ for $g\in L_{\infty }%
[ \mu ] ,$

\item the conditional entropy $S_{k,f}( \beta ) =\beta
E_{k,\beta f}[ f] -\ln Z_{k,\beta f},$

\item the conditional free energy $A_{k,f}( \beta ) =\beta
^{-1}\ln Z_{k,\beta f}$,

\item the conditional thermal variance $\sigma _{k,\beta f}^{2}[ g%
] =E_{k,\beta f}[ ( g-E_{k,\beta f}[ g]
) ^{2}] $ for $g\in L_{\infty }[ \mu ] $. As
$\beta \rightarrow 0,$ this becomes

\item the conditional variance $\sigma _{k}^{2}[ g]
=E_{k}[
( g-E_{k}[ g] ) ^{2}] $ for $g\in L_{\infty }%
[ \mu ] $.
\end{itemize}

If we fix all variables except $x_{k},$ then $E_{k}$ just becomes an
ordinary expectation, and it becomes evident that all the previously
established relations also hold for the corresponding conditional
quantities; in particular, the conclusions of Theorem \ref{Theorem
Entropy
fluctuation integral},
\[
S_{k,f}( \beta ) =\int_{0}^{\beta }\int_{t}^{\beta }\sigma
_{k,sf}^{2}[ f]\,\mathrm{d}s\,\mathrm{d}t,
\]
and of Lemma \ref{Lemma variational entropy bound},
\[
S_{k,f}( \beta ) \leq E_{k,\beta f}\bigl[ \psi \bigl( -\beta
( f-f_{k}) \bigr) \bigr]\qquad \mbox{if }f_{k}\in
\mathcal{A}_{k}.
\]
Other members of $\mathcal{A}_{k}$ that will play a role in the sequel
are:

\begin{itemize}
\item the conditional supremum $( \sup_{k}g) ( \mathbf{x}%
) =\sup_{y\in \Omega _{k}}g( \mathbf{x}_{y,k}) $ for
$g\in L_{\infty }[ \mu ] $,

\item the conditional infimum $( \inf_{k}g) ( \mathbf{x}%
) =\inf_{y\in \Omega _{k}}g( \mathbf{x}_{y,k}) $ for
$g\in L_{\infty }[ \mu ] $ and

\item the conditional range ran$_{k}( g)
=\sup_{k}g-\inf_{k}g$ for $g\in L_{\infty }[ \mu ] .$
\end{itemize}

\subsection{Tensorization of entropy}

In the non-interacting case, when the energy function $f$ is a sum,
$f=\sum f_{k},$ with $f_{k}\in L_{\infty }[ \mu _{k}] $, it
is easily verified that $S_{k,f}( \beta ) (
\mathbf{x}) =S_{k,f}( \beta ) $ is independent of
$\mathbf{x}$ and that
%
\begin{equation}\label{Additive_entropy}
S_{f}( \beta ) =\sum_{k=1}^{n}S_{k,f}( \beta ).
\end{equation}
Equality no longer holds in the interacting, nonlinear case, but there
is a subadditivity property that is sufficient for the purpose of
concentration inequalities.

The tensorization inequality states that the total entropy is no
greater than the thermal average of the sum\vadjust{\goodbreak} of the conditional
entropies. In 1975, Elliott Lieb \cite{Lieb1975} gave a proof of this
result, which was probably known some time before, at least in the
classical setting relevant to our arguments.

\begin{lemma}\label{Lemma Convexity of KL divergence}
Let $h,g>0$ be bounded
measurable
functions on $\Omega $. Then, for any expectation~$E$,
\[
E[ h] \ln \frac{E[ h] }{E[ g] }\leq
E\biggl[ h\ln \frac{h}{g}\biggr].
\]
\end{lemma}

\begin{pf}
Define an expectation functional $E_{g}$ by $E_{g}[ h] =E[ gh%
] /E[ g] $. The function $\Phi ( t) =t\ln
t$ is convex for positive $t$, since $\Phi ^{\prime \prime }=1/t>0$.
Thus, by
Jensen's inequality,
\[
E[ h] \ln \frac{E[ h] }{E[ g] }=E[ g%
] \Phi \biggl( E_{g}\biggl[ \frac{h}{g}\biggr] \biggr) \leq E[ g%
] E_{g}\biggl[ \Phi \biggl( \frac{h}{g}\biggr) \biggr] =E\biggl[
h\ln \frac{h}{g}\biggr] .
\]
\upqed\end{pf}

\begin{theorem}\label{Theorem Entropy subadditivity}
\begin{equation}\label{Entropy_subadditivity}
S_{f}( \beta ) \leq E_{\beta f}\Biggl[
\sum_{k=1}^{n}S_{k,f}( \beta ) \Biggr].
\end{equation}
\end{theorem}

\begin{pf}
We denote the canonical density with $\rho $, so $\rho =\mathrm{e}^{\beta
f}/Z_{\beta f}$. Writing $\rho =\rho /E[ \rho ] $ as a~telescopic product and
using the previous lemma, we get
\begin{eqnarray*}
E\biggl[ \rho \ln \frac{\rho }{E[ \rho ] }\biggr] &=&E\Biggl[
\rho
\ln \prod_{k=1}^{n}\frac{E_{1}\cdots E_{k-1}[ \rho ] }{%
E_{1}\cdots E_{k-1}E_{k}[ \rho ] }\Biggr] \\
&=&\sum E\biggl[ E_{1}\cdots E_{k-1}[ \rho ] \ln \frac{E_{1}\cdots E_{k-1}%
[ \rho ] }{E_{1}\cdots E_{k-1}[ E_{k}[ \rho ]
]
}\biggr] \\
&\leq &\sum E\biggl[ \rho \ln \frac{\rho }{E_{k}[ \rho ]
}\biggr]
=E\Biggl[ \sum E_{k}\biggl[ \rho \ln \frac{\rho }{E_{k}[ \rho ] }%
\biggr] \Biggr] .
\end{eqnarray*}
From the definition of $\rho $, we then obtain
\begin{eqnarray*}
S_{f}( \beta ) &=&\beta E_{\beta f}[ f] -\ln
Z_{\beta f}=E\biggl[ \rho \ln \frac{\rho }{E[ \rho ] }\biggr]
\leq E\biggl[
\sum E_{k}\biggl[ \rho \ln \frac{\rho }{E_{k}[ \rho ] }\biggr] %
\biggr] \\
&=&E\biggl[ \sum_{k=1}^{n}\biggl( E_{k}\biggl[ \frac{\mathrm{e}^{\beta f}}{Z_{\beta f}}%
\ln \frac{\mathrm{e}^{\beta f}}{Z_{\beta f}}\biggr] -E_{k}\biggl[ \frac{\mathrm{e}^{\beta f}}{%
Z_{\beta f}}\biggr] \ln E_{k}\biggl[ \frac{\mathrm{e}^{\beta f}}{Z_{\beta
f}}\biggr]
\biggr) \biggr] \\
&=&Z_{\beta f}^{-1}\sum_{k=1}^{n}E[ E_{k}[ \mathrm{e}^{\beta f}]
S_{k,f}( \beta ) ] =Z_{\beta
f}^{-1}\sum_{k=1}^{n}E[ \mathrm{e}^{\beta f}S_{k,f}( \beta )
] \qquad \mbox{by (\ref{conditional_identity})} \\
&=&E_{\beta f}\Biggl[ \sum_{k=1}^{n}S_{k,f}( \beta ) \Biggr].
\end{eqnarray*}%
\upqed\end{pf}

\subsection{The Efron--Stein--Steele inequality}\vspace*{-1pt}

Combining (\ref{Entropy_subadditivity}) with Theorem \ref{Theorem
Entropy
fluctuation integral} and dividing by $\beta ^{2}$, we obtain
\[
\frac{1}{\beta ^{2}}\int_{0}^{\beta }\int_{t}^{\beta }\sigma
_{sf}^{2}[
f]\,\mathrm{d}s\,\mathrm{d}t\leq E_{\beta f}\Biggl[ \sum_{k=1}^{n}\frac{1}{\beta ^{2}}%
\int_{0}^{\beta }\int_{t}^{\beta }\sigma _{k,sf}^{2}[ f]\,\mathrm{d}s\,\mathrm{d}t%
\Biggr].
\]
Using the continuity properties of $\beta \mapsto \sigma _{\beta
f}^{2}[ f] $, which follow from Lemma \ref{General Simple
Lemma}, we can take
the limit as $\beta \rightarrow 0$ and multiply by $2$ to obtain
\[
\sigma ^{2}[ f] \leq E\Biggl[ \sum_{k}\sigma _{k}^{2}[ f%
] \Biggr] ,
\]
which is the well-known Efron--Stein--Steele inequality
\cite{Steele1986}.
Observe that we may drop the assumption $f\in L_{\infty }[ \mu ] $%
, but we still require the existence of exponential moments in an
interval containing zero, so the inequality so derived is slightly
weaker than the inequality in \cite{Steele1986}.\vspace*{-1pt}

\subsection{A modified logarithmic Sobolev inequality}\vspace*{-1pt}

Suppose we have a sequence of functions $f_{k}\in \mathcal{A}_{k}$, so that $%
f_{k}$ is independent of the $k$th coordinate. Combining (\ref{Entropy_subadditivity}) with Lemma \ref{Lemma variational entropy bound} and
using the identity $E_{\beta f}E_{k,\beta f}=E_{\beta f}$, we obtain
%
\begin{equation}\label{ModLogSobInequ}
S_{f}( \beta ) \leq E_{\beta f}\Biggl[ \sum_{k=1}^{n}\psi
\bigl( -\beta ( f-f_{k}) \bigr) \Biggr] ,
\end{equation}
which is the modified logarithmic Sobolev inequality proposed by
Massart \cite{Massart2000a,Massart2000}. Many consequences of this powerful
inequality have been explored (e.g., \cite{Massart2000a,Boucheron2003,Bousquet2001,Maurer2006,Boucheron2009}).
Here we will concentrate on
the consequences of combining the tensorization inequality with the
fluctuation representation of entropy in Theorem \ref{Theorem Entropy
fluctuation integral}. Since the fluctuation representation is an
identity, this combination is stronger than (\ref{ModLogSobInequ}) and
leads
to some results that apparently cannot be recovered from (\ref%
{ModLogSobInequ}). We will also re-derive some results that can be
derived from (\ref{ModLogSobInequ}) in cases where we believe that the
proposed method gives some additional insight.\vspace*{-1pt}

\subsection{Conditional thermal variance and exponential concentration}\vspace*{-1pt}

Theorems \ref{Theorem basic Herbst argument},  \ref{Theorem Entropy
subadditivity} and \ref{Theorem Entropy fluctuation integral} (applied
to the conditional entropy) form the backbone of the proposed method.
Combining them, we obtain the following generic concentration result:

\begin{theorem}
\label{Theorem general concentration}For any $\beta >0$ we have the
entropy
bound
\[
S_{f}( \beta ) \leq E_{\beta f}\Biggl[ \sum_{k=1}^{n}\int_{0}^{%
\beta }\int_{t}^{\beta }\sigma _{k,sf}^{2}[
f]\,\mathrm{d}s\,\mathrm{d}t\Biggr],\vadjust{\goodbreak}
\]
the bound on the log-Laplace transform
\[
\ln E\bigl[ \mathrm{e}^{\beta ( f-Ef) }\bigr] =\beta \int_{0}^{\beta }%
\frac{S_{f}( \gamma ) }{\gamma ^{2}}\,\mathrm{d}\gamma
\]
and the concentration inequality
\[
\Pr \{ f-Ef>t\} \leq \exp \biggl( \beta \int_{0}^{\beta }\frac{%
S_{f}( \gamma ) }{\gamma ^{2}}\,\mathrm{d}\gamma -\beta t\biggr) .
\]
\end{theorem}

The obvious strategy is to start by bounding the conditional thermal
variance $\sigma _{k,sf}^{2}[ f] $. Typically, this leads to
considerable simplifications and we will follow this method in the
sequel.

\section{Two classical concentration inequalities}

We begin with the derivation of two classical results: the bounded
difference inequality and a~similar result, which reduces to the
familiar Bennett inequality when $f$ is the sum of its arguments. These
inequalities are not new, but they are very useful. We obtain them in
their strongest forms and they provide a good illustration of our
proposed method.

For $a,b\in
\mathbb{R}
$, $a<b$ define $\zeta _{a,b}\dvtx%
\mathbb{R}
\rightarrow
\mathbb{R}
$ by
\[
\zeta _{a,b}( t) =( b-t) ( t-a).
\]%
We state some elementary facts without proof.

\begin{lemma}
\label{Lemma_elementary}\textup{(i)} If $X$ is a random variable with values in $%
[ a,b], $ then
\[
\sigma ^{2}[ X] \leq ( b-EX) ( EX-a)
=\zeta _{a,b}( EX) \leq \frac{( b-a) ^{2}}{4}.
\]

\textup{(ii)} The function $\zeta _{a,b}$ is non-increasing in $[(
a+b) /2,\infty )$.
\end{lemma}

\subsection{The bounded difference inequality}

By Lemma \ref{Lemma_elementary}(i), we get for all $s\in
\mathbb{R}
$ that $\sigma _{k,sf}^{2}[ f] \leq \mathrm{ran}_{k}^{2}(
f) /4,$ so by the first conclusion of Theorem \ref{Theorem
general concentration},
%
\begin{equation}\label{Sumofsquaredrangesbound}
S_{f}( \gamma ) \leq \frac{1}{4}\int_{0}^{\gamma
}\int_{t}^{\gamma }E_{\gamma f}\Biggl[
\sum_{k=1}^{n}\mathrm{ran}_{k}^{2}( f) \Biggr]\,\mathrm{d}s\,\mathrm{d}t\leq
\frac{\gamma ^{2}}{8}E_{\gamma f}[ R^{2}( f) ],
\end{equation}
where we introduced the abbreviation $R^{2}( f) :=\sum_{k=1}^{n}$%
ran$_{k}^{2}( f) $. Bounding the thermal expectation by the
uniform norm, we obtain from the third conclusion of Theorem
\ref{Theorem
general concentration} that for all $\beta >0$%
\[
\Pr \{ f-Ef>t\} \leq \exp \biggl( \beta \int_{0}^{\beta }\frac{%
S_{f}( \gamma )\,\mathrm{d}\gamma }{\gamma ^{2}}-\beta t\biggr) \leq
\exp \biggl( \frac{\beta ^{2}}{8}\Vert R^{2}( f)
\Vert _{\infty }-\beta t\biggr) .
\]%
Substitution of the minimizing value $\beta =4t/\Vert R^{2}(
f) \Vert _{\infty }$ gives
\[
\Pr \{ f-Ef>t\} \leq \exp \biggl( \frac{-2t^{2}}{\Vert
R^{2}( f) \Vert _{\infty }}\biggr),
\]%
which is the well-known bounded difference inequality (with correct
exponent) in the strong version given by McDiarmid
\cite{McDiarmid1998}, Theorem 3.7, where the supremum is outside of the
sum of squared
conditional ranges. Note that the result is vacuous for $f\notin L_{\infty }%
[ \mu ] $.

\subsection{A Bennett--Bernstein concentration inequality}

The proof of the bounded difference inequality relied on bounding the
thermal variance~$\sigma _{k,\beta f}( f) $ uniformly in $\beta $%
, using constraints on the conditional range of $f$. We now consider
the
case where we only use one constraint on the ranges, say $f-E_{k}[ f%
] \leq 1$, but we use information on the conditional variances.
This leads to a Bennett-type inequality as in~\cite{McDiarmid1998},
Theorem 3.8. To state it, we abbreviate the sum of conditional
variances of $f$ as
\[
\Sigma ^{2}( f) =\sum \sigma _{k}^{2}( f) .
\]

Again, we start with a bound on the thermal variance.

\begin{lemma}
\label{Lemma Bennet variance bound}Assume $f-Ef\leq 1$. Then, for $\beta >0,$%
\[
\sigma _{\beta f}^{2}( f) \leq \mathrm{e}^{\beta}\sigma ^{2}(
f).
\]
\end{lemma}

\begin{pf}
\begin{eqnarray*}
\sigma _{\beta f}^{2}( f) &=&\sigma _{\beta (
f-Ef) }^{2}( f-Ef) =E_{\beta ( f-Ef)}[ (
f-Ef) ^{2}] -\bigl( E_{\beta ( f-Ef) }[ f-Ef] \bigr) ^{2} \\
&\leq &E_{\beta ( f-Ef) }[ ( f-Ef) ^{2}] =%
\frac{E[ ( f-Ef) ^{2}\mathrm{e}^{\beta ( f-Ef) }] }{E%
[ \mathrm{e}^{\beta ( f-Ef) }] } \\
&\leq &E\bigl[ ( f-Ef) ^{2}\mathrm{e}^{\beta ( f-Ef)
}\bigr]\qquad\mbox{use Jensen on denominator} \\
&\leq &\mathrm{e}^{\beta }E[ ( f-Ef) ^{2}] \qquad\mbox{use hypothesis.%
}
\end{eqnarray*}%
\upqed\end{pf}

Next we bound the total entropy $S_{f}( \beta ) .$

\begin{lemma}
\label{Lemma Bennett variance sum}Assume that $f-E_{k}f\leq 1$ for all
$k\in
\{ 1,\ldots,n\} $. Then, for $\beta >0,$
\[
S_{f}( \beta ) \leq ( \beta \mathrm{e}^{\beta }-\mathrm{e}^{\beta
}+1)E_{\beta f}[ \Sigma ^{2}( f) ] .
\]
\end{lemma}

\begin{pf}
Using the first conclusion of Theorem \ref{Theorem general
concentration}
and the previous lemma, we get
\[
S_{f}( \beta ) \leq E_{\beta f}\Biggl[ \sum_{k=1}^{n}\int_{0}^{%
\beta }\!\!\!\int_{t}^{\beta }\sigma _{k,sf}^{2}[ f]\,\mathrm{d}s\,\mathrm{d}t\Biggr]
\leq \int_{0}^{\beta }\!\!\!\int_{t}^{\beta }\mathrm{e}^{s}\,\mathrm{d}s\,\mathrm{d}tE_{\beta f}[
\Sigma ^{2}( f) ].\
\]
The conclusion follows from the elementary formula
\[
\int_{0}^{\beta }\!\!\!\int_{t}^{\beta }\mathrm{e}^{s}\,\mathrm{d}s\,\mathrm{d}t=\int_{0}^{\beta }(
\mathrm{e}^{\beta }-\mathrm{e}^{t})\,\mathrm{d}t=\beta \mathrm{e}^{\beta }-\mathrm{e}^{\beta }+1.
\]
\upqed\end{pf}

Now we can prove our version of Bennett's inequality.

\begin{theorem}\label{Theorem Bernstein inequality}
Assume $f-E_{k}f\leq 1,\forall k$. Let $%
t>0$ and denote $V=\Vert \Sigma ^{2}( f) \Vert
_{\infty }$. Then
\begin{eqnarray*}
\Pr \{ f-E[ f] >t\} &\leq &\exp \bigl( -V\bigl(
(1+tV^{-1}) \ln ( 1+tV^{-1}) -tV^{-1}\bigr) \bigr) \\
&\leq &\exp \biggl( \frac{-t^{2}}{2V+2t/3}\biggr) .
\end{eqnarray*}
\end{theorem}

\begin{pf}
Fix $\beta >0$. Recall the definition of the function $\psi $ in (\ref%
{Definition of psi}) and observe that
\[
\int_{0}^{\beta }\frac{\gamma \mathrm{e}^{\gamma }-\mathrm{e}^{\gamma }+1}{\gamma
^{2}}\,\mathrm{d}\gamma =\beta ^{-1}( \mathrm{e}^{\beta }-\beta -1) =\beta
^{-1}\psi ( \beta ),
\]
because $(\mathrm{d}/\mathrm{d}\gamma ) ( \gamma ^{-1}( \mathrm{e}^{\gamma
}-1) ) =\gamma ^{-2}( \gamma \mathrm{e}^{\gamma }-\mathrm{e}^{\gamma
}+1) $ and $\lim_{\gamma \rightarrow 0}\gamma ^{-1}(
\mathrm{e}^{\gamma
}-1) =1$. Theorem~\ref{Theorem general concentration} and Lemma \ref%
{Lemma Bennett variance sum} combined with a uniform bound then give
\begin{eqnarray*}
\ln E\mathrm{e}^{\beta ( f-Ef) } &=&\beta \int_{0}^{\beta }\frac{%
S_{f}( \gamma )\,\mathrm{d}\gamma }{\gamma ^{2}} \\
&\leq &\beta \biggl( \int_{0}^{\beta }\frac{\gamma \mathrm{e}^{\gamma }-\mathrm{e}^{\gamma }+1}{%
\gamma ^{2}}\,\mathrm{d}\gamma \biggr) \Vert \Sigma ^{2}( f)
\Vert _{\infty }=\psi ( \beta ) V.
\end{eqnarray*}
So, by Markov's inequality, we have for any $\beta >0$ that $\Pr \{ f-E%
[ f] >t\} \leq \exp ( \psi ( \beta )
V-\beta t) $. Substitution of $\beta =\ln ( 1+tV^{-1})
$ gives the first inequality; the second is Lemma 2.4 in
\cite{McDiarmid1998}.
\end{pf}

Observe that $f$ is assumed bounded above by the hypotheses of the
theorem. The existence of exponential moments $E[ \mathrm{e}^{\beta
f}] $ is needed
only for $\beta \geq 0$, so the assumption $f\in L_{\infty }[ \mu %
] $ can be dropped in this case.

\section{Exploiting monotonicity}

Sometimes an appropriately chosen bound on the conditional thermal variance $%
\sigma _{k,sf}^{2}[ f] $\vspace*{2pt} can be shown to have a monotonicity
property in the variable $s$, which can be used to find a bound uniform
in the the region of integration. The remaining part of the fluctuation
integral then just becomes $\beta ^{2}/2$, which leads to sub-Gaussian
tail estimates, just as for the bounded difference inequality. In this
section, we give three examples.

\subsection{Functions with large conditional expectations}

The following is our first novel result, the proof of which is hardly
more difficult than that of the bounded difference inequality. It
depends on the assumption that the conditional expectations are
consistently in the upper
halves of the conditional ranges for all $k$ and all configurations $%
x_{1},\ldots,x_{k-1},x_{k+1},\ldots,x_{n}$ of the conditioning data. If this
condition is met, the result is much stronger than the bounded
difference inequality, and, for large deviations $t$, also much
stronger than Bennett's inequality.

\begin{theorem}
\label{Theorem_Coherent}Suppose that
%
\begin{equation} \label{Coherence_condition}
E_{k}[ f] \geq \frac{\sup_{k}f+\inf_{k}f}{2}\qquad\forall k\in
\{ 1,\ldots,n\}
\end{equation}
and let
\[
A=\biggl\Vert \sum_{k}\Bigl( \sup_{k}f-E_{k}[ f] \Bigr)
\Bigl( E_{k}[ f] -\inf_{k}f\Bigr) \biggr\Vert _{\infty }.
\]%
Then%
\[
\Pr \{ f-Ef>t\} \leq \mathrm{e}^{-t^{2}/(2A)}.
\]
\end{theorem}

\begin{pf}
By Lemma \ref{General Simple Lemma}, the function $\beta \mapsto E_{k,\beta f}%
[ f] $ is non-decreasing, so for $\beta \geq 0$ we have
\[
E_{k}[ f] \in \biggl[ \frac{\sup_{k}f+\inf_{k}f}{2},E_{k,\beta f}%
[ f] \biggr] .
\]%
Since the function $\zeta _{\inf_{k}f,\sup_{k}f}$ (of Lemma \ref{Lemma_elementary}) is non-increasing in this interval, we get
\[
\sigma _{k,\beta f}^{2}( f) \leq \zeta
_{\inf_{k}f,\sup_{k}f}( E_{k,\beta f}[ f] ) \leq
\zeta _{\inf_{k}f,\sup_{k}f}( E_{k}[ f] ) ,
\]
and, from the first conclusion of Theorem \ref{Theorem general concentration}%
,%
\begin{eqnarray*}
S_{f}( \beta ) &\leq &E_{\beta f}\Biggl[ \sum_{k=1}^{n}\int_{0}^{%
\beta }\int_{t}^{\beta }\sigma _{k,sf}^{2}[ f]\,\mathrm{d}s\,\mathrm{d}t\Biggr]
\leq \frac{\beta ^{2}}{2}E_{\beta f}\Biggl[ \sum_{k=1}^{n}\zeta
_{\inf_{k}f,\sup_{k}f}( E_{k}[ f] ) \Biggr] \\
&\leq &\frac{\beta ^{2}}{2}\Biggl\Vert \sum_{k=1}^{n}\zeta
_{\inf_{k}f,\sup_{k}f}( E_{k}[ f] ) \Biggr\Vert
_{\infty }=\frac{\beta ^{2}A}{2}.
\end{eqnarray*}
The result now follows as in the proof of the bounded difference
inequality.
\end{pf}

\subsection{Monotonicity of variational bounds on the thermal variance}

A related strategy first finds a simple variational bound on the
conditional
thermal variance. We have%
\[
\sigma ^{2}[ g] =\min_{t\in
\mathbb{R}
}E[ ( g-t) ^{2}] \leq E[ ( g-c) ^{2}%
]\qquad\forall c\in
\mathbb{R}
.
\]%
Applied to the conditional thermal variance, this translates to%
\begin{equation} \label{Trivial variance bound}
\sigma _{k,\beta f}^{2}[ f] \leq E_{k,\beta f}[ (
f-f_{k}) ^{2}]\qquad\forall f_{k}\in
\mathcal{A}_{k}.
\end{equation}

We will use $\inf_{k}f$ for $f_{k}$ and combine this observation with
the following.

\begin{proposition}\label{Proposition_monotonicity}
The function $\beta \mapsto E_{k,\beta f}%
[ ( f-\inf_{k}f) ^{2}] $ is non-decreasing.
\end{proposition}

\begin{pf}
Write $h=f-\inf_{k}f$ and define a real function $\xi $ by $\xi (
t) =( \max \{ t,0\} ) ^{2}$. Since $h\geq
0,$ we have
\[
E_{k,\beta f}\Bigl[ \Bigl( f-\inf_{k}f\Bigr) ^{2}\Bigr] =E_{k,\beta
( f-\inf_{k}f) }\Bigl[ \Bigl( f-\inf_{k}f\Bigr) ^{2}\Bigr]
=E_{k,\beta h}[ \xi ( h) ] .
\]
By Lemma \ref{General Simple Lemma}, we obtain
\[
\frac{\mathrm{d}}{\mathrm{d}\beta }E_{\beta h}[ \xi ( h) ] =E_{\beta h}%
[ \xi ( h) h] -E_{\beta h}[ \xi ( h) %
] E_{\beta h}[ h] \geq 0,
\]
where the last inequality uses the well-known fact that for any expectation $%
E[ \xi ( h) h] \geq E[ \xi ( h)
] E[ h] $ whenever $\xi $ is a non-decreasing function.
\end{pf}

A first consequence is a lower tail bound somewhat similar to
Bernstein's inequality, Theorem~\ref{Theorem Bernstein inequality}.

\begin{theorem}
\label{Theorem Positive 2002}Let $t>0$ and denote
\[
W=\biggl\Vert \sum_{k}E_{k}\Bigl( f-\inf_{k}f\Bigr) ^{2}\biggr\Vert
_{\infty }.
\]
Then
\[
\Pr \{ E[ f] -f>t\} \leq \exp \biggl( \frac{-t^{2}}{2W}%
\biggr) .
\]
\end{theorem}

\begin{pf}
We use inequality (\ref{Trivial variance bound}) and Proposition \ref{Proposition_monotonicity} to get for $s\geq 0$
\[
\sigma _{k,-sf}^{2}[ f] \leq E_{k,-sf}\Bigl[ \Bigl(
f-\inf_{k}f\Bigr) ^{2}\Bigr] \leq E_{k}\Bigl[ \Bigl(
f-\inf_{k}f\Bigr) ^{2}\Bigr] .
\]
We therefore obtain from Theorem \ref{Theorem general concentration}
\begin{eqnarray*}
S_{-f}( \beta ) &\leq &E_{-\beta f}\Biggl[ \sum_{k=1}^{n}%
\int_{0}^{\beta }\int_{t}^{\beta }\sigma _{k,-sf}^{2}[ f]\,\mathrm{d}s\,\mathrm{d}t%
\Biggr] \leq \frac{\beta ^{2}}{2}E_{-\beta f}\Biggl[ \sum_{k=1}^{n}E_{k}%
\Bigl( f-\inf_{k}f\Bigr) ^{2}\Biggr] \\
&\leq &\frac{\beta ^{2}W}{2}
\end{eqnarray*}
and then proceed as in the proof of the bounded difference inequality.
\end{pf}

If we take the function $f$ to be an average of real random variables,
then Theorem \ref{Theorem Positive 2002} reduces to an inequality given
in \cite{Pinelis1989} and \cite{Maurer2003}. In \cite{Maurer2003} it is
argued that for very heterogeneous variables this inequality is
superior to Bernstein's inequality. Similar arguments apply to the
present, more general case.

When we apply the same method to obtain upper tail bounds we arrive at
a surprisingly powerful result. To state it, we introduce worst-case
variance proxies, which will play an important role in the sequel.

\begin{definition}
Let $g\in L_{\infty }[ \mu ] $. The worst-case variance proxy of $%
g$ is the function $Dg\in L_{\infty }[ \mu ] $ defined by
\[
Dg=\sum_{k}\Bigl( g-\inf_{k}g\Bigr) ^{2}.
\]
\end{definition}

The function $Dg$ is a local measure of the sensitivity of $g$ to
modifications of its individual arguments.

\begin{lemma}
\label{Lemma Entropy bound Upper}We have, for $\beta >0$,
\[
S_{f}( \beta ) \leq \frac{\beta ^{2}}{2}E_{\beta f}[ Df%
] .
\]
\end{lemma}

\begin{pf}
We use inequality (\ref{Trivial variance bound}) and
Proposition \ref{Proposition_monotonicity} to get for $0\leq s\leq \beta $
\[
\sigma _{k,sf}^{2}[ f] \leq E_{k,sf}\Bigl[ \Bigl(
f-\inf_{k}f\Bigr) ^{2}\Bigr] \leq E_{k,\beta f}\Bigl[ \Bigl(
f-\inf_{k}f\Bigr) ^{2}\Bigr] .
\]
So, using Theorem \ref{Theorem general concentration} again,
\begin{eqnarray*}
S_{f}( \beta ) &\leq &E_{\beta f}\Biggl[ \sum_{k=1}^{n}\int_{0}^{%
\beta }\int_{t}^{\beta }\sigma _{k,sf}^{2}[ f]\,\mathrm{d}s\,\mathrm{d}t\Biggr]
\leq \frac{\beta ^{2}}{2}E_{\beta f}\Biggl[ \sum_{k=1}^{n}E_{k,\beta
f}\Bigl(
f-\inf_{k}f\Bigr) ^{2}\Biggr] \\
&=&\frac{\beta ^{2}}{2}E_{\beta f}\Biggl[ \sum_{k=1}^{n}\Bigl(
f-\inf_{k}f\Bigr) ^{2}\Biggr] ,
\end{eqnarray*}
where we used the identity $E_{\beta f}E_{k,\beta f}=E_{\beta f}$ in
the last equation.
\end{pf}

The usual arguments now immediately lead to the following.

\begin{theorem}
\label{Theorem upper tail}With $t>0,$
\[
\Pr \{ f-E[ f] >t\} \leq \exp \biggl( \frac{-t^{2}}{%
2\Vert Df\Vert _{\infty }}\biggr) .
\]
\end{theorem}

In \cite{Maurer2006} this result was derived from
inequality (\ref{ModLogSobInequ}), and it is shown that it improves the exponent on
upper tail bounds derived from Talagrand's convex distance inequality
in many cases, as for shortest travelling salesmen paths, Steiner trees
and the
eigenvalues of random symmetric matrices. Here we only give one example of how $%
\Vert Df\Vert _{\infty }$ may be bounded and consider a
convex Lipschitz function $f$ defined on the cube $[ 0,1]
^{n}$. For simplicity, we assume $f$ to be differentiable.

Let $x\in [ 0,1] ^{n}$ and suppose that for some fixed $k$
there is $y\in [ 0,1] $ such that $f(
\mathbf{x}_{y,k}) \leq f( \mathbf{x}) $. Then by
convexity (using really only the fact
that $f$ is separately convex in each coordinate),
\[
f( \mathbf{x}) -f( \mathbf{x}_{y,k}) \leq
\langle \mathbf{x}-\mathbf{x}_{y,k},\partial f(
\mathbf{x}) \rangle
_{%
\mathbb{R}
^{n}}=( x_{k}-y) \partial _{k}f( \mathbf{x}) \leq
\vert \partial _{k}f( \mathbf{x}) \vert .
\]
We therefore have $f( \mathbf{x}) -\inf_{y}f( \mathbf{x}%
_{y,k}) \leq \vert \partial _{k}f( \mathbf{x})
\vert $ and
\[
Df( \mathbf{x}) =\sum_{k=1}^{n}\Bigl( f(
\mathbf{x}) -\inf_{y}f( \mathbf{x}_{y,k}) \Bigr)
^{2}\leq \Vert
\partial f( \mathbf{x}) \Vert _{%
\mathbb{R}
^{n}}^{2}\leq \Vert f\Vert _{\mathrm{Lip}}^{2}.
\]
In combination with Theorem \ref{Theorem upper tail} we obtain upper
tail bounds for $f$ with an exponent twice as good as obtained from the
convex distance inequality \cite{Ledoux2001}, Corollary 4.10, or an
earlier application of the entropy method \cite{Ledoux2001}, Theorem
5.9.

For a corresponding lower tail bound, we have to use an estimate
similar to what was used in the proof of Bennett's inequality.

\begin{lemma}
\label{Lemma Entropy Bound Lower}If $f-\inf_{k}f\leq 1,\forall k,$ then for $%
\beta >0,$%
\[
S_{-f}( \beta ) \leq \psi ( \beta ) E_{-\beta
f}[ Df] ,
\]%
with $\psi $ defined as in (\ref{Definition of psi}).
\end{lemma}

\begin{pf}
Let $k\in \{ 1,\ldots,n\} $. We write $h_{k}:=f-\inf_{k}f$. Then $%
h_{k}\in [ 0,1] $ and for $s\leq \beta $
\[
E_{k,-sh_{k}}[ h_{k}^{2}] =\frac{E_{k}[
h_{k}^{2}\mathrm{e}^{-\beta h_{k}}\mathrm{e}^{( \beta -s) h_{k}}]
}{E_{k}[ \mathrm{e}^{-\beta h_{k}}\mathrm{e}^{( \beta -s) h_{k}}]
}\leq \mathrm{e}^{( \beta -s) }\frac{E_{k}[ h_{k}^{2}\mathrm{e}^{-\beta
h_{k}}] }{E_{k}[ \mathrm{e}^{-\beta h_{k}}] }=\mathrm{e}^{( \beta
-s) }E_{k,-\beta h_{k}}[ h_{k}^{2}] .
\]
We therefore have%
\begin{eqnarray*}
\int_{0}^{\beta }\!\!\!\int_{t}^{\beta }E_{k,-sf}[ h_{k}^{2}]\,\mathrm{d}s\,\mathrm{d}t
&=&\int_{0}^{\beta }\!\!\!\int_{t}^{\beta }E_{k,-sh_{k}}[
h_{k}^{2}]\,\mathrm{d}s\,\mathrm{d}t \\
&\leq &\biggl( \int_{0}^{\beta }\!\!\!\int_{t}^{\beta }\mathrm{e}^{\beta
-s}\,\mathrm{d}s\,\mathrm{d}t\biggr) E_{k,-\beta h_{k}}[ h_{k}^{2}] =\psi (
\beta ) E_{k,-\beta f}[ h_{k}^{2}] ,
\end{eqnarray*}
where we used the formula
\[
\int_{0}^{\beta }\!\!\!\int_{t}^{\beta }\mathrm{e}^{-s}\,\mathrm{d}s\,\mathrm{d}t=1-\mathrm{e}^{-\beta }-\beta
\mathrm{e}^{-\beta }.
\]
Thus, using Theorem \ref{Theorem general concentration} and the identity $%
E_{-\beta f}E_{k,-\beta f}=E_{-\beta f},$
\begin{eqnarray*}
S_{-f}( \beta ) &\leq &E_{-\beta f}\biggl[
\sum_{k}\int_{0}^{\beta }\!\!\!\int_{t}^{\beta }\sigma _{k,-sf}^{2}[
f]\,\mathrm{d}s\,\mathrm{d}t\biggr] \leq E_{-\beta f}\biggl[ \sum_{k}\int_{0}^{\beta
}\!\!\!\int_{t}^{\beta }E_{k,-sf}[
h_{k}^{2}]\,\mathrm{d}s\,\mathrm{d}t\biggr] \\
&\leq &\psi ( \beta ) E_{-\beta f}\biggl[ \sum_{k}E_{k,-\beta f}%
[ h_{k}^{2}] \biggr] =\psi ( \beta ) E_{-\beta
f}[ Df] .
\end{eqnarray*}
\upqed\end{pf}

Lemmas \ref{Lemma Entropy bound Upper} and \ref{Lemma Entropy Bound
Lower} together with Theorem \ref{Theorem basic Herbst argument} imply
the
inequalities
%
\begin{equation}\label{Main Lemma Beta Positive}
\ln E\bigl[ \mathrm{e}^{\beta ( f-E[ f] ) }\bigr] \leq \frac{%
\beta }{2}\int_{0}^{\beta }E_{\gamma f}[ Df]\,\mathrm{d}\gamma
\end{equation}
and, if $f-\inf_{k}f\leq 1$ for all $k$, then
%
\begin{equation}\label{Main Lemma Beta Negative}
\ln E\bigl[ \mathrm{e}^{\beta ( E[ f] -f) }\bigr] \leq \frac{%
\psi ( \beta ) }{\beta }\int_{0}^{\beta }E_{-\gamma f}[ Df%
]\,\mathrm{d}\gamma ,
\end{equation}
where in the last inequality we also used the fact that $\gamma \mapsto
\psi ( \gamma ) /\gamma ^{2}$ is non-decreasing. Bounding
the thermal expectation with the uniform norm and substitution of
$\beta =\ln ( 1+t\Vert Df\Vert _{\infty }^{-1}) $
gives the following lower tail bound that can also be found in
\cite{Maurer2006}.

\begin{theorem}
\label{Theorem lower tail}If $f-\inf_{k}f\leq 1$ for all $k$, then for $t>0,$
\begin{eqnarray*}
\Pr \{ Ef-f>t\} & \leq& \exp \biggl( -\Vert Df\Vert
_{\infty }\biggl( \biggl( 1+\frac{t}{\Vert Df\Vert _{\infty }}%
\biggr) \ln \biggl( 1+\frac{t}{\Vert Df\Vert _{\infty }}\biggr) -%
\frac{t}{\Vert Df\Vert _{\infty }}\biggr) \biggr) \\
& \leq& \exp \biggl( \frac{-t^{2}}{2\Vert Df\Vert _{\infty }+2t/3}%
\biggr).
\end{eqnarray*}
\end{theorem}

The two inequalities (\ref{Main Lemma Beta Positive}) and (\ref{Main
Lemma Beta Negative}) are the keys to obtaining concentration
inequalities in terms of the worst-case variance proxy $Df$. Both
results can also be deduced from Massart's inequality
(\ref{ModLogSobInequ}) as shown in \cite{Maurer2006}. We do not claim
that the derivations given above are per se superior. We presented them
because they follow the same principles as the proofs of the
other results given above (the bounded difference inequality and Theorems \ref%
{Theorem Bernstein inequality}, \ref{Theorem_Coherent} and %
\ref{Theorem Positive 2002}), which do not follow from inequality (\ref%
{ModLogSobInequ}).

\section{Self-boundedness and canonical decoupling}

We conclude by presenting two general principles to extend the utility
of the proposed method. All the above applications of Theorem
\ref{Theorem general concentration} involved a chain of inequalities of
the form
\[
S_{\epsilon f}( \gamma ) \leq E_{\epsilon \gamma f}\Biggl[
\sum_{k=1}^{n}\int_{0}^{\gamma }\int_{t}^{\gamma }\sigma _{k,\epsilon sf}^{2}%
[ f]\,\mathrm{d}s\,\mathrm{d}t\Biggr] \leq \xi ( \gamma ) E_{\epsilon
\gamma f}[ G( f) ] ,
\]
where $\epsilon =1$ for upper tail results and $\epsilon =-1$ for lower
tail results, $\xi $ is some non-negative real function and $G(
f) $ is some function on $\Omega $ derived from $f$. For the
bounded difference
inequality, for example, $\xi ( \gamma ) =\gamma ^{2}/8$ and $%
G=R^{2}( f) $; for the Bennett inequality $\xi ( \gamma
) =\gamma \mathrm{e}^{\gamma }-\mathrm{e}^{\gamma }+1$ and $G( f)
=\Sigma ^{2}( f) $; for Theorem \ref{Theorem upper tail} we
had $\xi ( \gamma ) =\gamma ^{2}/2$ and $G( f)
=Df$; while for
the corresponding lower tail bound, Theorem \ref{Theorem lower tail}, we had $%
\xi ( \gamma ) =\psi ( \gamma ) $ and also
$G( f) =Df$, etc. Theorem \ref{Theorem general
concentration} is then
invoked to conclude that
%
\begin{equation}\label{Uniform_bound}
\ln E\mathrm{e}^{\epsilon \beta ( f-Ef) }\leq \beta \int_{0}^{\beta }\frac{%
\xi ( \gamma ) }{\gamma ^{2}}E_{\epsilon \gamma f}[
G( f) ]\,\mathrm{d}\gamma \leq \beta \Vert G( f)
\Vert _{\infty }\int_{0}^{\beta }\frac{\xi ( \gamma )\,\mathrm{d}\gamma }{\gamma ^{2}}.
\end{equation}

Here the uniform estimate $E_{\epsilon \beta f}[ G( f) %
] \leq \Vert G( f) \Vert _{\infty }$, while
being very simple, is somewhat loose. We now sketch how it can
sometimes be avoided by exploiting special properties of the thermal
expectation.

\subsection{Self-boundedness}

The first possibility we consider is that the function $G(
f) $ can be bounded in terms of the function $f$ itself, a
property referred to
as \textit{self-boundedness} \cite{Boucheron2009}. For example, if  $%
G( f) \leq f,$ then $E_{\gamma f}[ G( f)
] \leq E_{\gamma f}[ f] =( \mathrm{d}/\mathrm{d}\gamma ) \ln
E[ \exp ( \gamma f) ] $, and if the function $\xi
$ has some reasonable behavior, then the first integral in
(\ref{Uniform_bound}) above can be bounded by partial integration or
even more easily. As an example, we
apply this idea in the setting of Theorems \ref{Theorem upper tail} and \ref%
{Theorem lower tail}.

\begin{theorem}
\label{Theorem_selfbound}Suppose that there are non-negative numbers
$a,b$
such that $Df\leq af+b$. Then, for $t>0,$ we have
\[
\Pr \{ f-E[ f] >t\} \leq \exp \biggl( \frac{-t^{2}}{%
2( aE[ f] +b+at/2) }\biggr) .
\]
If, in addition, $a\geq 1$ and $f-\inf_{k}f\leq 1,\forall k\in \{
1,\ldots,n\} $, then
\[
\Pr \{ E[ f] -f>t\} \leq \exp \biggl( \frac{-t^{2}}{%
2( aE[ f] +b) }\biggr) .
\]
\end{theorem}

\begin{pf}
We only prove the lower tail bound; for the upper tail we refer to
\cite{Maurer2006}. As for the lower tail, it follows from (\ref{Main
Lemma Beta
Negative}) and Lemma \ref{General Simple Lemma} that
\begin{eqnarray*}
\ln E\bigl[ \mathrm{e}^{\beta ( E[ f] -f) }\bigr] &\leq &\frac{%
a\psi ( \beta ) }{\beta }\int_{0}^{\beta }E_{-\gamma f}[ f%
]\,\mathrm{d}\gamma +b\psi ( \beta ) =\frac{-a\psi ( \beta
) }{\beta }\ln Z_{-\beta f}+b\psi ( \beta ) \\
&=&\frac{-a\psi ( \beta ) }{\beta }\ln E\bigl[ \mathrm{e}^{\beta ( E%
[ f] -f) }\bigr] +\psi ( \beta ) (aE[ f] +b) .
\end{eqnarray*}
Rearranging gives%
\[
\ln E\bigl[ \mathrm{e}^{\beta ( E[ f] -f) }\bigr] \leq \frac{%
\psi ( \beta ) }{1+a\beta ^{-1}\psi ( \beta )
}( aE[ f] +b) \leq \frac{\beta ^{2}(
aE[ f] +b) }{2},
\]%
where one verifies that for $\beta >0$ and $a\geq 1$ we have $\psi
( \beta ) ( 1+a\beta ^{-1}\psi ( \beta )
) ^{-1}\leq \beta ^{2}/2$. The usual analysis with Markov's
inequality and optimization in $\beta $ conclude the proof.
\end{pf}

Recently Boucheron \textit{et al.} \cite{Boucheron2009} have given a refined
version of this result, where the condition $a\geq 1$ is improved to
$a\geq 1/3$ for the lower tail. There they also show that Theorems
\ref{Theorem_selfbound} and \ref{Theorem upper tail} together suffice
to derive a version of the convex distance inequality that differs from
Talagrand's original result only in that it has an inferior exponent.

It must be stressed that the same method of proof can be used to yield
self-bounded versions of all concentration inequalities derived from Theorem %
\ref{Theorem general concentration}, such as the bounded-difference and
Bennett inequalities.\vspace*{-3pt}

\subsection{Decoupling}\vspace*{-3pt}

A second method to avoid the uniform bound on the thermal expectation
uses decoupling. Recall that for any two probability measures $\nu $
and $\mu $ and a measurable function $g$ we have
\[
E_{x\sim \nu }[ g( x) ] \leq KL( \nu ,\mu
) +\ln E_{x\sim \mu }\mathrm{e}^{g( x) },
\]
which can be regarded as an instance of convex duality and easily
verified directly from the definition of the Kullback--Leibler
divergence. Applying this inequality when $\nu $ is the canonical
ensemble and $\mu $ is the a priori measure, we obtain for any $\theta
>0$
\[
S_{\epsilon f}( \beta ) \leq \xi ( \beta ) \theta
^{-1}E_{\epsilon \beta f}[ \theta G( f) ] \leq
\xi ( \beta ) \theta ^{-1}\bigl( S_{\epsilon f}( \beta
) +\ln E[ \exp ( \theta G( f) )
] \bigr) .
\]
For values of $\beta $ and $\theta $ where $\theta >\xi ( \beta
)
$ we obtain
\[
S_{\epsilon f}( \beta ) \leq \frac{\xi ( \beta ) }{%
\theta -\xi ( \beta ) }\ln E[ \exp ( \theta
G( f) ) ] .
\]
Hence, if we can control the upwards deviations of $G( f) $
(or some suitable bound thereof), we obtain concentration inequalities
for $f$ in terms of the expectation of $G( f) $ (or the
bound thereof). Again, this method, which was proposed in
\cite{Boucheron2003}, can be applied to all the versions of $G(
f) $ we introduced above and combined with all methods to control
the upwards deviation of $G( f) $, which leads to a
proliferation of concentration inequalities. Perhaps not all of these
deserve to be documented. We just quote a corresponding result in
\cite{Maurer2010} that uses $G( f) =Df$ and combines with
self-boundedness.

\begin{theorem}
\label{Theorem_Concentration}Suppose that there is $g\in L_{\infty
}[
\mu ] $ and $a\geq 1$ such that $0\leq f\leq g$, $Df\leq ag$ and $%
Dg\leq ag$. Then, for $t>0$,
\[
\Pr \{ f-Ef>t\} \leq \exp \biggl( \frac{-t^{2}}{4aE[
g] +3at/2}\biggr) .\vadjust{\goodbreak}
\]
If, in addition, $f-\inf_{k}f\leq 1$ for all $k$, then
\[
\Pr \{ Ef-f>t\} \leq \exp \biggl( \frac{-t^{2}}{4aE[
g] +at}\biggr) .
\]
\end{theorem}

In \cite{Maurer2010} the theorem is used to show that the concentration
of eigenvalues ($f$) of the Gram matrix of a sample of independent,
bounded random vectors in a Hilbert space is controlled by the size of
the largest eigenvalue ($g$).

\section{A glossary of notation}

We conclude with a tabular summary of notation.\\[5pt]
\begin{tabular}{@{}ll@{}}
$\Omega =\prod_{k=1}^{n}\Omega _{k}$ & underlying (product-) probability space.\\[1pt]
$\mu =\bigotimes _{k=1}^{n}\mu _{k}$ & (product-) probability measure on $%
\Omega $. \\[1pt]
$X_{k}$ & random variable distributed as $\mu _{k}$ in $\Omega _{k}$. \\[1pt]
$f\in L_{\infty }[ \mu ] $ & fixed function (negative
energy)
under investigation. \\[1pt]
$g\in L_{\infty }[ \mu ] $ & generic function. \\[1pt]
$E[ g] =\int_{\Omega }g\,\mathrm{d}\mu $ & expectation of $g$ in $\mu $. \\[1pt]
$\sigma ^{2}[ g] =E[ ( g-E[ g] ) ^{2}%
] $ & variance of $g$ in $\mu $. \\[1pt]
$\beta =1/T$ & inverse temperature. \\[1pt]
$E_{\beta f}[ g] =E[ g\mathrm{e}^{\beta f}] /E[ \mathrm{e}^{\beta f}%
] $ & thermal expectation of $g$. \\[1pt]
$Z_{\beta f}=E[ \mathrm{e}^{\beta f}] $ & partition function. \\[1pt]
$S_{f}( \beta ) =\beta E_{\beta f}[ f] -\ln
Z_{\beta
f}.$ & canonical entropy. \\[1pt]
$A_{f}( \beta ) =\frac{1}{\beta }\ln Z_{\beta f}$ &
Helmholtz
free energy. \\
$\sigma _{\beta f}^{2}( g) =E_{\beta f}[ ( g-E_{\beta f}%
[ g] ) ^{2}] $ & thermal variance of $g$. \\[1pt]
$\psi ( t) =\mathrm{e}^{t}-t-1$ &  \\[1pt]
$\mathbf{x}_{y,k}$ & vector $\mathbf{x}\in \Omega $ with $x_{k}$
replaced by
$y\in \Omega _{k}$. \\[1pt]
$E_{k}[ g] ( \mathbf{x}) =\int_{\Omega
_{k}}g( \mathbf{x}_{y,k})\,\mathrm{d}\mu _{k}( y) $ &
conditional expectation.\\[1pt]
$\mathcal{A}_{k}\subset L_{\infty }[ \mu ] $ & functions
independent of $k$th variable. \\[1pt]
$Z_{k,\beta f}=E_{k}[ \mathrm{e}^{\beta f}] $ & conditional partition
function. \\[1pt]
$E_{k,\beta f}[ g] =Z_{k,\beta f}^{-1}E_{k}[ g\mathrm{e}^{\beta f}%
] $ & conditional thermal expectation. \\[1pt]
$S_{k,f}( \beta ) =\beta E_{k,\beta f}[ g] -\ln
Z_{k,\beta f}$ & conditional entropy. \\[1pt]
$\sigma _{k,\beta f}^{2}[ g] =E_{k,\beta f}[ (
g-E_{k,\beta f}[ g] ) ^{2}] $ & conditional
thermal
variance. \\[1pt]
$\sigma _{k}^{2}[ g] =E_{k}[ ( g-E_{k}[
g]
) ^{2}] $ & conditional variance. \\[1pt]
$( \sup_{k}g) ( \mathbf{x}) =\sup_{y\in \Omega
_{k}}g( \mathbf{x}_{y,k}) $ & conditional supremum. \\[1pt]
$( \inf_{k}g) ( \mathbf{x}) =\inf_{y\in \Omega
_{k}}g( \mathbf{x}_{y,k}) $ & conditional infimum. \\[1pt]
ran$_{k}( g) =\sup_{k}g-\inf_{k}g$ & conditional range. \\[1pt]
$R^{2}( g) =\sum_{k}$ran$_{k}^{2}( g) $ & sum of
conditional square ranges. \\[1pt]
$\Sigma ^{2}( g) =\sum_{k}\sigma _{k}^{2}[ g] $ &
sum
of conditional variances. \\[1pt]
$Dg=\sum_{k}( g-\inf_{k}g) ^{2}$ & worst case variance proxy.%
\end{tabular}


\printhistory

\end{document}